\documentclass[11pt,a4paper]{amsart}

\usepackage[margin=1in]{geometry}
\usepackage{amsmath,amssymb,amsthm}
\IfFileExists{mathtools.sty}{\usepackage{mathtools}}{}
\usepackage[authoryear,round]{natbib}
\usepackage{booktabs,graphicx,array}
\usepackage{xcolor}
\usepackage{microtype}
\usepackage[hidelinks]{hyperref}
\hypersetup{
  pdftitle={Sub-Gaussian HDLSS Principal Component Analysis},
  pdfauthor={Yugo Nakayama},
  pdfsubject={High-probability theory and numerical evidence for NRM, CDM, and block noise reduction},
  pdfkeywords={HDLSS, PCA, sub-Gaussian, NRM, CDM, bNRM}
}
\graphicspath{{figures/}}
\numberwithin{equation}{section}

\newtheorem{theorem}{Theorem}[section]
\newtheorem{lemma}[theorem]{Lemma}
\newtheorem{proposition}[theorem]{Proposition}
\newtheorem{corollary}[theorem]{Corollary}
\theoremstyle{definition}
\newtheorem{definition}[theorem]{Definition}
\newtheorem{assumption}[theorem]{Assumption}
\theoremstyle{remark}
\newtheorem{remark}[theorem]{Remark}
\newtheorem{example}[theorem]{Example}

\DeclareMathOperator{\tr}{tr}

\DeclareMathOperator{\diag}{diag}

\newcommand{\R}{\mathbb{R}}

\newcommand{\op}{\mathrm{op}}
\newcommand{\F}{\mathrm{F}}
\newcommand{\ind}{\mathbf{1}}
\newcommand{\eps}{\varepsilon}
\newcommand{\wh}[1]{\widehat{#1}}
\newcommand{\wt}[1]{\widetilde{#1}}
\newcommand{\norm}[1]{\left\lVert#1\right\rVert}
\newcommand{\abs}[1]{\left\lvert#1\right\rvert}

\newcommand{\authoraffiliation}{Independent researcher}
\makeatletter
\def\@setauthors{%
  \begingroup
  \trivlist
  \centering\footnotesize \@topsep30\p@\relax
  \advance\@topsep by -\baselineskip
  \item\relax
  \author@andify\authors
  \def\\{\protect\linebreak}%
  {\large\scshape\authors\par}%
  \vspace{4pt}
  {\normalfont\normalsize\itshape\authoraffiliation\par}%
  \endtrivlist
  \endgroup
}
\makeatother

\title[Sub-Gaussian HDLSS PCA]{Sub-Gaussian HDLSS Principal Component Analysis:\\
High-Probability Theory and Numerical Evidence for NRM, CDM, and Block Noise Reduction}
\author{Yugo Nakayama}
\date{9 November 2026}

\begin{document}
\begin{abstract}
We study principal component analysis when the dimension $d$ is much larger than the sample size $n$ under an independent-factor sub-Gaussian model.  The dual sample covariance is decomposed into a finite-rank score matrix, an isotropic bulk shift, and a high-probability remainder governed by the bulk trace $\tau_0$, Frobenius scale $q_0$, and largest noise eigenvalue $b_0=\norm{\Lambda_0}_{\op}$.  This decomposition separates fixed-score HDLSS limits from joint limits in which $n$ also grows.  It shows that sequential noise reduction generally targets a score oracle containing a later-spike subtraction, while an oracle block correction avoids that subtraction when the signal rank is known.  A parallel cross-data-matrix bound explains the benefit and cost of sample splitting.  Reproducible Gaussian-factor simulations with 1,000 replications per configuration confirm the score-oracle approximation, the multi-spike subtraction, the separation between eigenvalue and eigenvector scales, and the population transition when $n$ grows.
\end{abstract}
\maketitle

\section{Introduction}\label{sec:intro}
\enlargethispage{2pt}
Principal component analysis (PCA) estimates the eigensystem of a covariance matrix.  Its geometric and probabilistic foundations go back to \citet{Pearson1901} and \citet{Hotelling1933}; standard accounts and a modern review are given by \citet{Jolliffe2002} and \citet{JolliffeCadima2016}.  In the classical regime, the dimension is fixed and the sample size tends to infinity.

When dimension and sample size grow proportionally, the Mar\v{c}enko--Pastur law and spiked-covariance theory describe bulk spectra, extreme eigenvalues, and eigenvector phase transitions \citep{MarchenkoPastur1967,Johnstone2001,BaikEtAl2005,BaikSilverstein2006,Paul2007}.  The present paper instead emphasizes the high-dimension, low-sample-size (HDLSS) regime $d\gg n$, where the sample covariance has rank at most $n-1$.  Its nonzero eigenvalues are more naturally studied through the $n\times n$ dual Gram matrix than through the $d\times d$ covariance matrix; see \citet{HallMarronNeeman2005}, the extension of the geometric representation by \citet{AhnEtAl2007}, \citet{JungMarron2009}, and the surveys of \citet{AoshimaYata2017,AoshimaEtAl2018}.

\citet{YataAoshima2012} developed the noise-reduction methodology (NRM), which subtracts an estimate of the high-dimensional noise floor from each empirical eigenvalue.  The cross-data-matrix methodology (CDM) of \citet{YataAoshima2010} instead splits the sample and crosses independent halves, replacing the positive within-sample noise shift by a mean-zero cross-noise term.  Related HDLSS developments include non-Gaussian PCA consistency \citep{YataAoshima2009}, intrinsic-dimension estimation \citep{YataAoshima2010Intrinsic}, first-component inference \citep{IshiiEtAl2016}, and empirical-eigenstructure bias analysis \citep{WangFan2017}.  The power-spiked analysis of \citet{YataAoshima2013} places these procedures in a common asymptotic framework.  Sparsity offers a complementary route to direction recovery, as illustrated by \citet{JohnstoneLu2009}; here we study dense directions and explicit noise-floor corrections.

The phrase ``under a sub-Gaussian assumption'' is not by itself precise enough for a quadratic-form proof.  A tail bound for every linear functional does not automatically provide Hanson--Wright concentration for $x^\top A x$.  We consequently distinguish directional sub-Gaussianity from an independent-factor representation and use the latter in the main results.

The central point is informational.  Increasing $d$ creates many bulk coordinates, so their aggregate Gram contribution can become nearly deterministic.  It does not create additional observations of the low-dimensional signal scores.  When $n$ is fixed, the score covariance remains random; even perfect bulk removal need not consistently estimate a population spike.

\subsection{Contributions}
The paper makes five contributions.
\begin{enumerate}
\item It gives a nonasymptotic decomposition of the dual covariance.  The bulk error is expressed by $\tau_0=\tr(\Lambda_0)$, $q_0=\norm{\Lambda_0}_{\F}$, and $b_0=\norm{\Lambda_0}_{\op}$.
\item It identifies the exact finite-score target of sequential NRM.  The target for component $i$ contains $-(r-i)^{-1}\sum_{j=i+1}^m s_j$, which need not vanish at fixed $n$.
\item Conditional on the correct signal rank $m$, it analyzes block NRM (bNRM), which estimates the common bulk floor from only the lower empirical eigenvalues.
\item It derives a comparable high-probability CDM bound and separates the cancellation of the mean shift from the variability introduced by sample splitting.
\item It reports deterministic, reproducible numerical experiments that distinguish population error from finite-score-oracle error rather than conflating the two targets.
\end{enumerate}

The results exchange distributional generality for explicit probabilities.  The independent-factor sub-Gaussian assumption is stronger than the fourth-moment conditions used in parts of the original NRM/CDM literature.  The claims below concern theorem form, proof transparency, and interpretation; they do not assert that every correction displayed here is absent from all previous covariance-denoising work.

\section{Model, regimes, and notation}\label{sec:setup}
\subsection{Dual representation}
Let $x_1,\ldots,x_n\in\R^d$ be independent, mean-zero observations with covariance $\Sigma_d$.  Define
\begin{equation}
X=[x_1,\ldots,x_n]\in\R^{d\times n},\qquad P_n=I_n-\frac1n\ind_n\ind_n^\top,\qquad r=n-1,
\end{equation}
and $Y=XP_n$.  The sample covariance and dual sample covariance are
\begin{equation}
S=\frac1rYY^\top,\qquad S_D=\frac1rY^\top Y=\frac1rP_nX^\top XP_n.
\end{equation}
If $Y=U\Gamma V^\top$ is a singular value decomposition, then $S=r^{-1}U\Gamma^2U^\top$ and $S_D=r^{-1}V\Gamma^2V^\top$; hence they have the same nonzero eigenvalues.

Fix $m<r$ and decompose
\begin{equation}\label{eq:covsplit}
\Sigma_d=H_s\Lambda_sH_s^\top+H_0\Lambda_0H_0^\top,
\end{equation}
where $H_s^\top H_0=0$, $\Lambda_s=\diag(\lambda_1,\ldots,\lambda_m)$, and $\Lambda_0=\diag(\lambda_{m+1},\ldots,\lambda_d)$.  Set
\begin{equation}\label{eq:bulksummaries}
\tau_0=\tr(\Lambda_0),\qquad q_0=\norm{\Lambda_0}_{\F}=\Bigl(\sum_{j=m+1}^d\lambda_j^2\Bigr)^{1/2},\qquad b_0=\norm{\Lambda_0}_{\op}=\max_{j>m}\lambda_j.
\end{equation}
Thus $b_0$ is explicitly the largest noise eigenvalue, equivalently the operator norm of the diagonal bulk matrix.  Two useful effective ranks are $r_{\mathrm{eff}}(\Lambda_0)=\tau_0/b_0$ and $r_2(\Lambda_0)=\tau_0^2/q_0^2$.

\subsection{Sub-Gaussian formulations}
For a real random variable $Z$, write
\begin{equation}
\norm{Z}_{\psi_2}=\inf\{t>0:\mathbb{E}\exp(Z^2/t^2)\le2\}.
\end{equation}
\begin{definition}[Directional sub-Gaussianity]
A mean-zero random vector $x$ with covariance $\Sigma$ is directionally $K$-sub-Gaussian if $\norm{u^\top x}_{\psi_2}\le K\sqrt{u^\top\Sigma u}$ for every $u\in\R^d$.
\end{definition}
\begin{assumption}[Independent-factor sub-Gaussian model, SG-I]\label{ass:sgi}
For each observation,
\begin{equation}\label{eq:factor}
x_k=H_s\Lambda_s^{1/2}z_{s,k}+H_0\Lambda_0^{1/2}z_{0,k}.
\end{equation}
All coordinates of $z_{s,k}$ and $z_{0,k}$ are mutually independent, have mean zero and variance one, and satisfy $\norm{z_{jk}}_{\psi_2}\le K$.  Observations are independent across $k$.
\end{assumption}
Assumption~\ref{ass:sgi} implies directional sub-Gaussianity up to an absolute constant, but the converse is false in general.  The distinction is essential because the proofs apply Hanson--Wright concentration to $z^\top\Lambda_0z$.

\subsection{Two HDLSS regimes}
We keep $m$ fixed and separate:
\begin{description}
\item[Fixed-score HDLSS (F-HDLSS)] $n$ and $m$ are fixed while $d\to\infty$.
\item[Joint HDLSS (J-HDLSS)] $d\to\infty$, $n=n_d\to\infty$, $n_d/d\to0$, and $m$ is fixed.
\end{description}
In F-HDLSS, the $m\times n$ signal-score matrix has a fixed number of columns and does not concentrate at its population covariance.  In J-HDLSS it can.

\subsection{Notation table}
Table~\ref{tab:notation} collects the symbols used throughout the paper.
\begin{table}[ht]
\centering
\caption{Principal notation.}\label{tab:notation}
\small
\begin{tabular}{@{}p{0.19\textwidth}p{0.72\textwidth}@{}}
\toprule
Symbol & Meaning \\
\midrule
$d,n,r$ & feature dimension, sample size, and centered rank $r=n-1$ \\
$P_n$ & centering projection $I_n-n^{-1}\ind_n\ind_n^\top$ \\
$X,Y$ & uncentered and centered data matrices; $Y=XP_n$ \\
$S,S_D$ & feature-space and dual sample covariance matrices \\
$m$ & signal rank, assumed fixed with $m<r$ \\
$H_s,H_0$ & orthonormal signal and bulk eigenvector matrices \\
$\Lambda_s,\Lambda_0$ & signal and bulk eigenvalue matrices \\
$\tau_0$ & bulk trace $\tr(\Lambda_0)$; the mean dual-noise shift scale \\
$q_0$ & bulk Frobenius norm $\norm{\Lambda_0}_{\F}$; quadratic fluctuation scale \\
$b_0$ & largest noise eigenvalue $\max_{j>m}\lambda_j=\norm{\Lambda_0}_{\op}$ \\
$Z_s,Z_0$ & standardized signal and bulk factor-score matrices \\
$A_s$ & full-sample signal score oracle $r^{-1}P_nZ_s^\top\Lambda_sZ_sP_n$ \\
$s_i$ & nonzero eigenvalues of $A_s$ \\
$a_i$ & ordered nonzero empirical eigenvalues of $S_D$ \\
$\wt\lambda_i,\wh\lambda_i^B$ & sequential NRM and oracle-rank bNRM estimators \\
$C_D,C_s$ & cross-data matrix and its signal score-oracle part \\
$\Delta_N,\Delta_C$ & NRM and CDM high-probability remainder bounds \\
\bottomrule
\end{tabular}
\end{table}

\section{Concentration and perturbation tools}\label{sec:tools}
\begin{lemma}[Hanson--Wright inequality]\label{lem:hw}
Let $z=(z_1,\ldots,z_p)^\top$ have independent, mean-zero coordinates with $\norm{z_j}_{\psi_2}\le K$.  For every symmetric $A$ and $t\ge0$,
\begin{equation}
\Pr\!\left(\abs{z^\top Az-\mathbb{E}z^\top Az}>t\right)
\le2\exp\!\left[-c\min\!\left\{\frac{t^2}{K^4\norm{A}_{\F}^2},\frac{t}{K^2\norm{A}_{\op}}\right\}\right].
\end{equation}
\end{lemma}
This is the standard Hanson--Wright inequality; see \citet{RudelsonVershynin2013} and \citet{HsuKakadeZhang2012}.

\begin{lemma}[Bilinear form bound]\label{lem:bilinear}
Let $z,w\in\R^p$ be independent and satisfy the conditions of Lemma~\ref{lem:hw}.  For every symmetric $A$ and $u\ge0$, with probability at least $1-2e^{-u}$,
\begin{equation}\label{eq:bilinear}
\abs{z^\top Aw}\le CK^2\left(\norm{A}_{\F}\sqrt{u}+\norm{A}_{\op}u\right).
\end{equation}
\end{lemma}
\begin{proof}
Polarization gives
\begin{equation}
4z^\top Aw=(z+w)^\top A(z+w)-(z-w)^\top A(z-w).
\end{equation}
The coordinates of $z\pm w$ are independent across their coordinate index and have $\psi_2$ norm at most $2K$.  Apply Lemma~\ref{lem:hw} to both quadratic forms, use a union bound, and absorb numerical constants into $C$.
\end{proof}

\begin{lemma}[Weyl inequalities]\label{lem:weyl}
For symmetric $A,B\in\R^{p\times p}$ with ordered eigenvalues $a_j,b_j$,
\begin{equation}
\max_j\abs{a_j-b_j}\le\norm{A-B}_{\op}.
\end{equation}
For rectangular matrices of the same size, $\abs{\sigma_j(A)-\sigma_j(B)}\le\norm{A-B}_{\op}$ for every singular-value index $j$.
\end{lemma}
These perturbation facts, together with Davis--Kahan bounds for eigenspaces, are standard; see \citet{StewartSun1990} and \citet{YuWangSamworth2015}.

\section{High-probability dual decomposition}\label{sec:dual}
Let $Z_s=[z_{s,1},\ldots,z_{s,n}]$ and $Z_0=[z_{0,1},\ldots,z_{0,n}]$.  Orthogonality in \eqref{eq:covsplit} implies the exact identity
\begin{equation}\label{eq:gramsplit}
X^\top X=Z_s^\top\Lambda_sZ_s+Z_0^\top\Lambda_0Z_0.
\end{equation}
\begin{theorem}[Dual bulk decomposition]\label{thm:dual}
Assume SG-I and let $u\ge0$.  Put
\begin{equation}
L_N(u,n)=u+\log(2n^2),\qquad
\Delta_N(u)=\frac{CK^2n}{r}\left\{q_0\sqrt{L_N(u,n)}+b_0L_N(u,n)\right\}.
\end{equation}
Then, with probability at least $1-2e^{-u}$,
\begin{equation}\label{eq:dualdecomp}
S_D=A_s+\frac{\tau_0}{r}P_n+E_N,
\qquad A_s=\frac1rP_nZ_s^\top\Lambda_sZ_sP_n,
\qquad \norm{E_N}_{\op}\le\Delta_N(u).
\end{equation}
\end{theorem}
\begin{proof}
Write $G_0=Z_0^\top\Lambda_0Z_0$.  For its diagonal entries, coordinate independence and unit variance give $\mathbb{E}(G_0)_{jj}=\tau_0$; observation independence gives $\mathbb{E}(G_0)_{jk}=0$ when $j\ne k$.  Hence $\mathbb{E}G_0=\tau_0I_n$.

Lemma~\ref{lem:hw} bounds every diagonal entry of $G_0-\tau_0I_n$, and Lemma~\ref{lem:bilinear} bounds every off-diagonal entry.  With $L=L_N(u,n)$, a union bound over at most $n^2$ entries gives, after changing the absolute constant,
\begin{equation}
\max_{j,k}\abs{(G_0-\tau_0I_n)_{jk}}
\le CK^2\{q_0\sqrt L+b_0L\}
\end{equation}
with probability at least $1-2e^{-u}$.  Since $\norm{M}_{\op}\le n\max_{j,k}\abs{M_{jk}}$ and $P_n$ is an orthogonal projection,
\begin{equation}
\norm{r^{-1}P_n(G_0-\tau_0I_n)P_n}_{\op}\le\Delta_N(u).
\end{equation}
Substitution into $S_D=r^{-1}P_nX^\top XP_n$ using \eqref{eq:gramsplit} and $P_n^2=P_n$ proves the result.
\end{proof}

\begin{corollary}[Fixed-$n$ spherical bulk]\label{cor:spherical}
In F-HDLSS, suppose $q_0/\tau_0\to0$ and $b_0/\tau_0\to0$.  Then
\begin{equation}
\frac{r}{\tau_0}(S_D-A_s)\xrightarrow{p}P_n.
\end{equation}
\end{corollary}
\begin{proof}
Fix $\eps>0$ and choose $u$ so that $2e^{-u}<\eps/2$.  Since $n$ is fixed, $L_N(u,n)$ is fixed and $\Delta_N(u)/\tau_0\to0$.  The event in Theorem~\ref{thm:dual} then implies the stated convergence, while its complement has probability below $\eps/2$.
\end{proof}

\begin{remark}
The factor $n$ in $\Delta_N$ comes from the transparent but conservative maximum-entry inequality.  Direct matrix concentration may improve the dependence on $n$; this is not needed for the fixed-$n$ conclusions.
\end{remark}

\section{NRM, its finite-score target, and block correction}\label{sec:nrm}
Work on the $r$-dimensional space $\ind_n^\perp$, where $P_n$ acts as the identity.  Let $a_1\ge\cdots\ge a_r$ be the nonzero eigenvalues of $S_D$.  Sequential NRM is
\begin{equation}\label{eq:nrm}
\wt\lambda_i=a_i-\frac1{r-i}\sum_{j=i+1}^r a_j,\qquad 1\le i\le r-1.
\end{equation}
Let $s_1\ge\cdots\ge s_m\ge0$ be the nonzero eigenvalues of $A_s$ and set $s_{m+1}=\cdots=s_r=0$.  These $s_i$ are random finite-score oracle eigenvalues.

\begin{definition}[Oracle block noise reduction]\label{def:bnrm}
If $m$ is known, define
\begin{equation}\label{eq:bnrm}
\wh\nu_m=\frac1{r-m}\sum_{j=m+1}^r a_j,
\qquad \wh\lambda_i^B=a_i-\wh\nu_m,\quad i=1,\ldots,m.
\end{equation}
When $m$ is estimated, all results below are conditional on the event $\wh m=m$.
\end{definition}

\begin{theorem}[Finite-score oracle bounds]\label{thm:oracle}
Suppose $\norm{E_N}_{\op}\le\Delta$ in \eqref{eq:dualdecomp}.  Then
\begin{equation}\label{eq:bnrmbound}
\max_{1\le i\le m}\abs{\wh\lambda_i^B-s_i}\le2\Delta,
\end{equation}
and, for sequential NRM,
\begin{equation}\label{eq:nrmtarget}
\abs{\wt\lambda_i-\left(s_i-\frac1{r-i}\sum_{j=i+1}^m s_j\right)}\le2\Delta,
\qquad 1\le i\le m.
\end{equation}
\end{theorem}
\begin{proof}
On $\ind_n^\perp$, compare $S_D$ with $B=A_s+(\tau_0/r)I_r$.  The eigenvalues of $B$ are $b_j=s_j+\tau_0/r$ for $j\le m$ and $b_j=\tau_0/r$ for $j>m$.  Lemma~\ref{lem:weyl} gives $\max_j\abs{a_j-b_j}\le\Delta$.

The lower-eigenvalue average of $B$ is exactly $\tau_0/r$.  Averaging cannot increase the maximum perturbation, so $\abs{\wh\nu_m-\tau_0/r}\le\Delta$, and \eqref{eq:bnrmbound} follows by the triangle inequality.

For sequential NRM, direct substitution gives
\begin{equation}
 b_i-\frac1{r-i}\sum_{j=i+1}^r b_j
 =s_i-\frac1{r-i}\sum_{j=i+1}^m s_j.
\end{equation}
The NRM functional is $2$-Lipschitz in the maximum norm, which proves \eqref{eq:nrmtarget}.
\end{proof}

\begin{corollary}[Comparable spikes at fixed $n$]\label{cor:twospike}
Suppose $m=2$, $s_1=s_2=s>0$, and $\Delta/s\xrightarrow{p}0$.  Then
\begin{equation}
\frac{\wt\lambda_1}{s}\xrightarrow{p}1-\frac1{r-1}=\frac{r-2}{r-1},
\qquad
\frac{\wh\lambda_1^B}{s}\xrightarrow{p}1.
\end{equation}
\end{corollary}
The corollary does not say that sequential NRM fails in every multi-spike regime.  It says that later-spike subtraction is not automatically negligible when $n$ is fixed.

\section{From score oracles to population eigenvalues}\label{sec:population}
The nonzero eigenvalues $s_i$ of $A_s$ equal those of
\begin{equation}\label{eq:Ms}
M_s=\Lambda_s^{1/2}\left(\frac1rZ_sP_nZ_s^\top\right)\Lambda_s^{1/2}.
\end{equation}
This reduction does not make $s_i$ equal to $\lambda_i$ when $n$ is fixed.

\begin{example}[Rank-one fixed-score counterexample]\label{ex:rankone}
With no bulk noise, let $x_k=\sqrt\lambda z_kh$, where $z_k\stackrel{\mathrm{iid}}{\sim}N(0,1)$ and $\norm h=1$.  The only nonzero empirical eigenvalue is
\begin{equation}
\wh\lambda=\frac\lambda r\sum_{k=1}^n(z_k-\bar z)^2,
\qquad \frac{\wh\lambda}{\lambda}\sim\frac{\chi_r^2}{r}.
\end{equation}
The distribution does not involve $d$ and therefore cannot converge to one as $d\to\infty$ at fixed $n$.
\end{example}
Thus perfect bulk removal can recover a random score oracle but not, in general, a population eigenvalue.

\begin{lemma}[Signal-score covariance concentration]\label{lem:scoreconcentration}
Under SG-I, for every $u\ge0$, with probability at least $1-2e^{-u}$,
\begin{equation}\label{eq:scoreconcentration}
\norm{\frac1rZ_sP_nZ_s^\top-I_m}_{\op}
\le\eps_s(u):=CK^2\left(\sqrt{\frac{m+u}{r}}+\frac{m+u}{r}\right).
\end{equation}
\end{lemma}
This is the fixed-dimensional sub-Gaussian sample-covariance bound; see \citet{Vershynin2018}.  Effective-rank refinements for sample covariance operators are developed by \citet{KoltchinskiiLounici2017}.  Centering changes $n$ to $r=n-1$.

\begin{theorem}[Population eigenvalue bound for bNRM]\label{thm:population}
On the common event in Theorem~\ref{thm:dual} and Lemma~\ref{lem:scoreconcentration},
\begin{equation}\label{eq:populationbound}
\abs{\wh\lambda_i^B-\lambda_i}\le2\Delta_N(u)+\lambda_1\eps_s(u),
\qquad i=1,\ldots,m.
\end{equation}
Consequently, in J-HDLSS, if $u_n\to\infty$, $u_n=o(n)$,
\begin{equation}
\frac{\Delta_N(u_n)}{\lambda_i}\to0,
\qquad
\frac{\lambda_1}{\lambda_i}\left(\sqrt{\frac{m+u_n}{r}}+\frac{m+u_n}{r}\right)\to0,
\end{equation}
then $\wh\lambda_i^B/\lambda_i\xrightarrow{p}1$.
\end{theorem}
\begin{proof}
Equation~\eqref{eq:Ms} and submultiplicativity yield $\norm{M_s-\Lambda_s}_{\op}\le\lambda_1\eps_s(u)$.  Weyl's inequality gives $\abs{s_i-\lambda_i}\le\lambda_1\eps_s(u)$.  Combine this with Theorem~\ref{thm:oracle} and use the triangle inequality.  The asymptotic conclusion follows after division by $\lambda_i$; the failure probability tends to zero because $u_n\to\infty$.
\end{proof}

Sequential NRM requires the additional population-scale condition
\begin{equation}
\frac1{(r-i)\lambda_i}\sum_{j=i+1}^m\lambda_j\to0.
\end{equation}
For fixed $m$, $n\to\infty$, and bounded spike ratios, this term is $O(n^{-1})$.

\section{Eigenvectors require the trace scale}\label{sec:eigenvectors}
Eigenvalue recovery depends on fluctuations of quadratic forms and is governed by $q_0$ and $b_0$.  Feature-space directions contain the full bulk energy and therefore depend on $\tau_0$.

Consider $m=1$, with signal eigenvector $h$ and eigenvalue $\lambda$.  Let $a=P_nz_s$, $v=a/\norm a$ when $a\ne0$, and define the oracle reconstruction $q=XP_nv=Xv$.  Then
\begin{equation}
q=\sqrt\lambda\,\norm a\,h+\eta,
\qquad \eta=H_0\Lambda_0^{1/2}Z_0v,
\qquad h^\top\eta=0.
\end{equation}
\begin{proposition}[Exact rank-one angle identity]\label{prop:angle}
Let $\wh h=q/\norm q$.  Then
\begin{equation}\label{eq:angle}
\sin^2\angle(\wh h,h)=\frac{\norm\eta^2}{\lambda\norm a^2+\norm\eta^2}.
\end{equation}
Conditionally on $v$, $\mathbb{E}(\norm\eta^2\mid v)=\tau_0$ and, with probability at least $1-2e^{-u}$,
\begin{equation}\label{eq:anglenoise}
\abs{\norm\eta^2-\tau_0}\le CK^2(q_0\sqrt u+b_0u).
\end{equation}
\end{proposition}
\begin{proof}
Orthogonality gives $\norm q^2=\lambda\norm a^2+\norm\eta^2$ and $\abs{h^\top\wh h}^2=\lambda\norm a^2/\norm q^2$, proving \eqref{eq:angle}.  Conditional on $v$, each coordinate of $Z_0v$ is mean-zero, variance-one, and sub-Gaussian, while coordinates are independent across the bulk index.  Hence $\norm\eta^2=(Z_0v)^\top\Lambda_0(Z_0v)$, and Lemma~\ref{lem:hw} gives \eqref{eq:anglenoise}.
\end{proof}

At fixed $n$, a sufficient scale for the oracle angle to vanish is roughly $\tau_0/\lambda\to0$, whereas eigenvalue separation may hold under $q_0/\lambda\to0$.  For an isotropic bulk, $\tau_0\asymp d$ but $q_0\asymp\sqrt d$.  Transferring Proposition~\ref{prop:angle} from the oracle score direction to an empirical eigenvector additionally requires an eigengap and a Davis--Kahan argument \citep{DavisKahan1970,YuWangSamworth2015}.

\section{Cross-data matrices and comparison}\label{sec:cdm}
Split the observations into independent groups of sizes $n_1,n_2$, set $r_\ell=n_\ell-1$, and let $P_\ell$ be the corresponding centering matrices.  Define
\begin{equation}\label{eq:cdm}
C_D=\frac1{\sqrt{r_1r_2}}P_1X_1^\top X_2P_2=C_s+E_C,
\end{equation}
where
\begin{equation}
C_s=\frac1{\sqrt{r_1r_2}}P_1Z_{s,1}^\top\Lambda_sZ_{s,2}P_2,
\qquad
E_C=\frac1{\sqrt{r_1r_2}}P_1Z_{0,1}^\top\Lambda_0Z_{0,2}P_2.
\end{equation}
There is no $\tau_0P/r$ term because the two groups are independent.

\begin{theorem}[CDM cross-noise bound]\label{thm:cdm}
Let $L_C(u)=u+\log(2n_1n_2)$ and
\begin{equation}
\Delta_C(u)=CK^2\sqrt{\frac{n_1n_2}{r_1r_2}}
\{q_0\sqrt{L_C(u)}+b_0L_C(u)\}.
\end{equation}
Under SG-I, with probability at least $1-2e^{-u}$,
\begin{equation}
\norm{E_C}_{\op}\le\Delta_C(u),
\qquad
\abs{\sigma_i(C_D)-\sigma_i(C_s)}\le\Delta_C(u).
\end{equation}
\end{theorem}
\begin{proof}
Every uncentered cross entry is $z_{0,j;1}^\top\Lambda_0z_{0,k;2}$, has mean zero, and satisfies Lemma~\ref{lem:bilinear}.  A union bound over $n_1n_2$ entries makes the bound simultaneous.  The inequality $\norm B_{\op}\le\norm B_{\F}\le\sqrt{n_1n_2}\max_{j,k}\abs{B_{jk}}$, contraction by the centering matrices, and rectangular Weyl perturbation complete the proof.
\end{proof}

\begin{proposition}[From the split score oracle to population spikes]\label{prop:cdmpop}
Set $A_\ell=r_\ell^{-1/2}Z_{s,\ell}P_\ell$ and $G_\ell=A_\ell A_\ell^\top$.  If $\norm{G_\ell-I_m}_{\op}\le\eps_\ell<1$, then
\begin{equation}
\abs{\sigma_i(C_s)-\lambda_i}\le C\lambda_1(\eps_1+\eps_2).
\end{equation}
\end{proposition}
\begin{proof}
Write $U_\ell=G_\ell^{-1/2}A_\ell$, so $U_\ell U_\ell^\top=I_m$ and
$C_s=U_1^\top G_1^{1/2}\Lambda_sG_2^{1/2}U_2$.  Partial isometries preserve its nonzero singular values.  Functional calculus gives $\norm{G_\ell^{1/2}-I_m}_{\op}\le\eps_\ell$, up to an absolute factor on $[0,1)$; expansion of $G_1^{1/2}\Lambda_sG_2^{1/2}-\Lambda_s$ and singular-value perturbation prove the claim.
\end{proof}

\begin{table}[ht]
\centering
\caption{Theoretical contrast under SG-I.}\label{tab:nrmcdm}
\small
\begin{tabular}{@{}p{0.24\textwidth}p{0.34\textwidth}p{0.34\textwidth}@{}}
\toprule
Feature & NRM / bNRM & CDM \\
\midrule
Mean bulk effect & Positive shift $\tau_0/r$, estimated and removed & Zero mean cross-noise \\
Random bulk scale & $q_0$ and $b_0$ & Same $q_0$ and $b_0$ scales \\
Observations & One full-sample Gram matrix & Split score samples \\
Multiple spikes & Sequential subtraction; bNRM avoids it if $m$ is correct & Joint singular values, no sequential term \\
Fixed-$n$ target & Full-sample score oracle & Split score oracle \\
Population consistency & Requires $n\to\infty$ in general & Requires both split score covariances to concentrate \\
\bottomrule
\end{tabular}
\end{table}

\section{Failure modes and power-spiked interpretation}\label{sec:boundaries}
\subsection{Directional sub-Gaussianity is insufficient}
Let $u$ be uniform on the sphere of radius $\sqrt d$ in $\R^d$ and let $R$ be independent with $R^2=1/2$ or $3/2$, each with probability $1/2$.  Then $x=Ru$ is isotropic and directionally sub-Gaussian, but $\norm x^2/d=R^2$ never concentrates at one.  Thus directional linear tails alone do not imply the bulk Gram concentration in Theorem~\ref{thm:dual}.

\subsection{Rank misspecification and repeated eigenvalues}
If $\wh m<m$, an omitted spike enters the lower-eigenvalue average and causes over-subtraction.  If $\wh m>m$, some noise eigenvalues are retained as signal.  The bNRM bound is therefore conditional on correct rank information.  When population eigenvalues repeat, individual eigenvectors are not identifiable; the appropriate target is the spectral projection onto the repeated eigenspace.

\subsection{Heavy tails}
Finite fourth moments do not imply sub-Gaussian tails.  Under heavy tails, the exponential inequalities here may fail, and robust Gram entries, truncation, median-of-means constructions, or robust covariance operators are needed.  This limitation distinguishes SG-I from parts of the Aoshima--Yata theory that use fourth-moment assumptions.

\subsection{Power-spiked models}
Suppose $\lambda_i=a_id^{\alpha_i}$ for $i\le m$ and $0<c_-\le\lambda_j\le c_+$ for $j>m$.  Then $\tau_0\asymp d$, $q_0\asymp\sqrt d$, and $b_0\asymp1$.  Ignoring logarithmic factors, the bulk fluctuation suggests eigenvalue recovery when $\alpha_i>1/2$, while fixed-$n$ rank-one direction recovery requires approximately $\alpha_i>1$.  These are heuristic phase boundaries: boundary exponents require sharper constants or limiting distributions, and growing $n$ relaxes the direction scale to approximately $\tau_0/(n\lambda_i)$.

\section{Numerical experiments}\label{sec:numerical}
\subsection{Design and reproducibility}
We ran four preregistered-style experiments under the Gaussian special case of SG-I with an isotropic unit bulk.  Every configuration used 1,000 Monte Carlo replications, master seed 20261109, and independent child streams created by NumPy's \texttt{SeedSequence}.  No draw was discarded and no tuning rule was selected after inspecting the results.  The complete script, summary CSV, seed manifest, and vector PDF figures are included in \texttt{code/}, \texttt{results/}, and \texttt{figures/}.  The largest Monte Carlo standard error among the reported mean errors was $0.0098$.

For an estimate $\wh\lambda_i$, we report population relative error $\abs{\wh\lambda_i-\lambda_i}/\lambda_i$ and score-oracle relative error $\abs{\wh\lambda_i-s_i}/s_i$.  For CDM, $s_i$ is replaced by the corresponding split score oracle $\sigma_i(C_s)$.  Unless stated otherwise, $n=20$.

\subsection{Fixed-score eigenvalue experiment}
\begin{figure}[ht]
\centering
\includegraphics[width=0.96\textwidth]{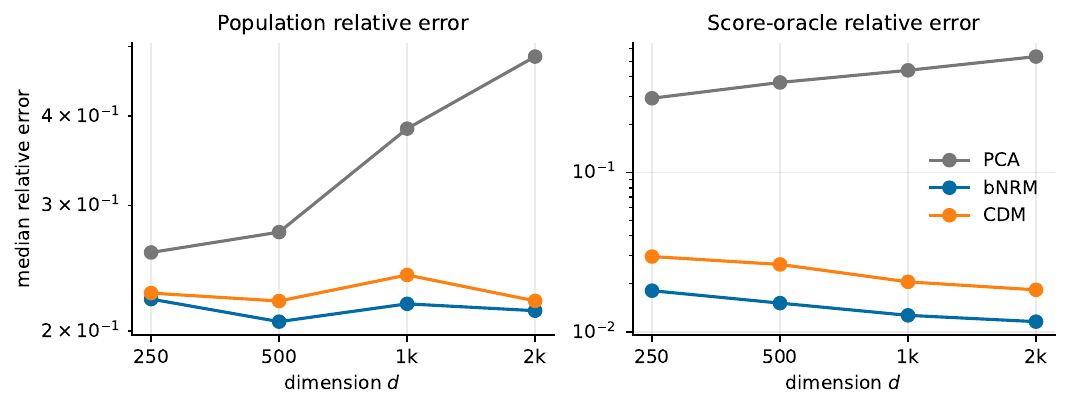}
\caption{Fixed-score experiment.  Medians over 1,000 replications.  bNRM and CDM approach their respective score oracles, whereas their population errors plateau because $n$ is fixed.}\label{fig:fixedscore}
\end{figure}

The first experiment used $m=1$, $\lambda=d^{0.7}$, and $d\in\{250,500,1000,2000\}$.  Ordinary PCA retains the bulk shift: its median population relative error rose from $0.258$ to $0.484$.  bNRM's score-oracle error fell from $0.0181$ to $0.0116$, and CDM's fell from $0.0296$ to $0.0183$.  In contrast, their population errors did not approach zero: the bNRM median was $0.222$ at $d=250$ and $0.214$ at $d=2000$; the corresponding CDM values were $0.226$ and $0.220$.  Figure~\ref{fig:fixedscore} therefore visualizes the central distinction: increasing $d$ improves bulk removal relative to the random score oracle but does not create new signal-score observations.

\subsection{Sequential subtraction with two spikes}
\begin{figure}[ht]
\centering
\includegraphics[width=0.62\textwidth]{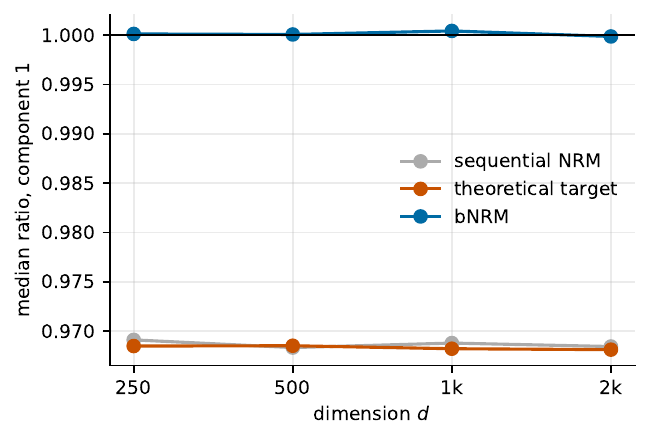}
\caption{Two equal population spikes.  Sequential NRM follows its later-spike-subtracted score target, while bNRM remains centered at the first score-oracle eigenvalue.}\label{fig:multispike}
\end{figure}

The second experiment used $m=2$, equal population spikes $\lambda_1=\lambda_2=d^{0.9}$, and the same dimensions.  For component one, Figure~\ref{fig:multispike} compares $\wt\lambda_1/s_1$ with the sample-specific theoretical target $1-s_2/\{(r-1)s_1\}$ and with $\wh\lambda_1^B/s_1$.  At $d=2000$, their medians were $0.96844$, $0.96813$, and $0.99987$, respectively.  The median absolute discrepancy between sequential NRM and its target decreased from $0.00435$ to $0.00201$, almost exactly matching the bNRM score-oracle discrepancy.  The sequential limit is not $17/18$ here because the two finite-sample score eigenvalues are ordered and typically unequal even when the population spikes coincide; Theorem~\ref{thm:oracle} predicts the realized ratio $s_2/s_1$.

\subsection{Eigenvalue and direction scales}
\begin{figure}[ht]
\centering
\includegraphics[width=0.62\textwidth]{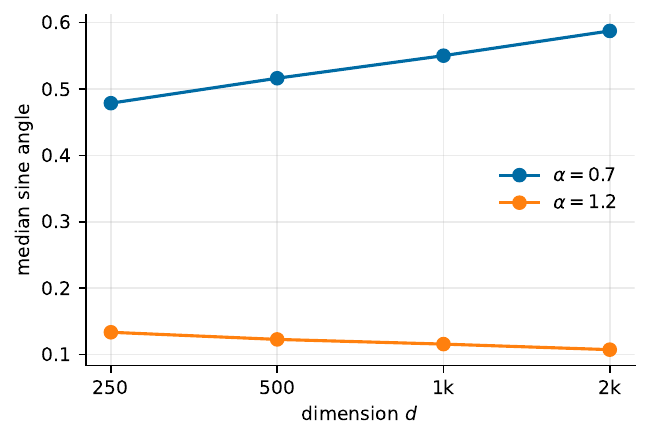}
\caption{Leading empirical direction under an isotropic bulk.  The weak spike is increasingly dominated by total bulk energy; the strong spike improves with dimension.}\label{fig:directions}
\end{figure}

The third experiment used one spike $\lambda=d^\alpha$.  With $\alpha=0.7$, the median sine angle of the leading empirical direction increased from $0.479$ to $0.587$ as $d$ grew from 250 to 2000.  With $\alpha=1.2$, it decreased from $0.134$ to $0.107$.  Thus $\alpha=0.7>1/2$ is sufficient for improving eigenvalue-oracle accuracy in the first experiment but not for direction recovery, whereas $\alpha=1.2>1$ overcomes the trace scale.  Figure~\ref{fig:directions} agrees with the phase separation in Proposition~\ref{prop:angle}.

\subsection{Joint HDLSS transition}
Finally, we set $n=\sqrt d$ at $d\in\{400,900,1600,2500\}$, with $m=1$ and $\lambda=d^{0.7}$.  The bNRM population-error median decreased from $0.221$ to $0.137$, while its score-oracle-error median decreased from $0.0162$ to $0.00437$.  The former remains larger because it contains score-covariance variability, but unlike the fixed-$n$ experiment it decreases as the number of score observations grows.

\begin{figure}[ht]
\centering
\includegraphics[width=0.62\textwidth]{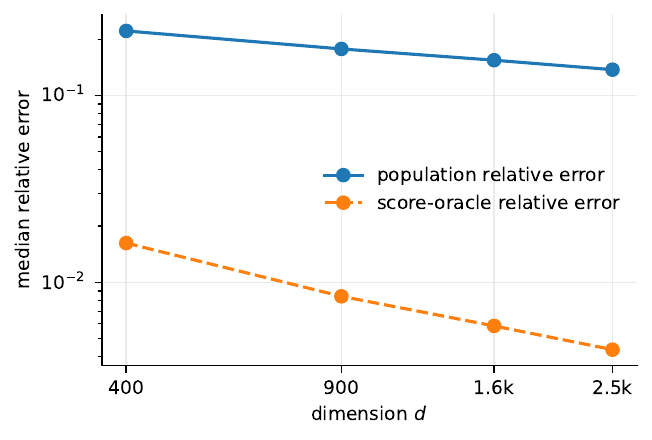}
\caption{J-HDLSS sequence with $n=\sqrt d$.  Both error components decrease, but score-covariance variability remains the dominant finite-sample contribution.}\label{fig:joint}
\end{figure}

\begin{table}[ht]
\centering
\caption{Selected numerical results (medians over 1,000 replications).}\label{tab:numerical}
\small
\begin{tabular}{@{}llrr@{}}
\toprule
Experiment & Metric & Initial $d$ & Final $d$ \\
\midrule
Fixed $n$, bNRM & population relative error & 0.222 & 0.214 \\
Fixed $n$, bNRM & score-oracle relative error & 0.0181 & 0.0116 \\
Fixed $n$, CDM & population relative error & 0.226 & 0.220 \\
Fixed $n$, CDM & score-oracle relative error & 0.0296 & 0.0183 \\
Two spikes, sequential NRM & estimate / score oracle & 0.9691 & 0.9684 \\
Two spikes, bNRM & estimate / score oracle & 1.0001 & 0.9999 \\
Direction, $\alpha=0.7$ & sine angle & 0.479 & 0.587 \\
Direction, $\alpha=1.2$ & sine angle & 0.134 & 0.107 \\
J-HDLSS, bNRM & population relative error & 0.221 & 0.137 \\
J-HDLSS, bNRM & score-oracle relative error & 0.0162 & 0.00437 \\
\bottomrule
\end{tabular}
\end{table}

\section{Discussion}\label{sec:discussion}
The gain from SG-I is a finite-$(d,n)$ failure probability.  The bulk contribution is not hidden inside an $o_p(1)$ term: its mean is $\tau_0$, its quadratic fluctuation is $q_0$, and its largest-coordinate obstruction is $b_0$.  The cost is equally important.  SG-I excludes many dependent designs, including spatial, image, genomic, and time-series dependence.  Extending the proof requires a quadratic-form concentration property, not merely a sub-Gaussian label.

Relative to representative NRM/CDM results, the present formulation changes four aspects.  It displays explicit remainders, replaces a required power law by trace/Frobenius/operator summaries, states fixed-$n$ targets as random score oracles, and makes sequential multi-spike subtraction visible.  These changes clarify interpretation rather than establishing a universal priority claim.

Important open problems are: a rank-selection oracle inequality for bNRM; dependent Hanson--Wright conditions; repeated-split CDM with dependence between splits; robust sub-Weibull or finite-moment variants; matrix-level bounds sharper in $n$ than the entrywise argument, potentially using matrix concentration tools such as those surveyed by \citet{Tropp2015}; and minimax lower bounds formalizing the fixed-score information limit.

\section{Conclusion}\label{sec:conclusion}
Under independent sub-Gaussian factors, the dual covariance separates into a finite-rank score term, an isotropic bulk shift, and a controlled remainder.  Sequential NRM and bNRM estimate the shift from one Gram matrix; CDM prevents the shift from appearing by crossing independent sample halves.  Their bulk errors depend on the same Frobenius and operator scales, but their score errors differ because CDM splits the sample.

Fixed-$n$ and growing-$n$ limits must not be conflated.  At fixed $n$, the natural targets are random score oracles.  Population eigenvalue consistency becomes available when the score covariance also concentrates.  Sequential NRM has an additional later-spike subtraction, while oracle-rank bNRM avoids it at the price of knowing the signal rank.  Finally, directions require a stronger signal than eigenvalues because they must overcome the full bulk energy $\tau_0$, not only the quadratic fluctuation $q_0$.  The numerical results reproduce all four distinctions with recorded seeds and code included in the source package.

\section*{Acknowledgments}
Parts of the drafting and technical editing of this manuscript were assisted by large language models, including Anthropic’s Claude Opus 5.5. The author takes full responsibility for the correctness of all mathematical statements, proofs, and numerical results.

\bibliographystyle{plainnat}
\bibliography{references}
\end{document}